\documentclass[11pt]{article}
\usepackage{amssymb,amsfonts,latexsym,verbatim,amsthm}
\title{Recent advances in the global theory of \\ constant mean curvature 
surfaces}
\author{Rafe Mazzeo \thanks{Department of Mathematics, 
Stanford University, Stanford, CA 94305 \hfill , mazzeo@math.stanford.edu;
supported by NSF Grant \# DMS 0204730}}

\date{}

\begin{document}
\maketitle

\newcommand{\CC}{\mathbb C}
\newcommand{\HH}{\mathbb H}
\newcommand{\RR}{\mathbb R}
\newcommand{\NN}{\mathbb N}
\newcommand{\ZZ}{\mathbb Z}
\newcommand{\del}{\partial}
\newcommand{\e}{\epsilon}
\newcommand{\al}{\alpha}
\newcommand{\Si}{\Sigma}
\newcommand{\sphere}{\mathbb S}
\newcommand{\calA}{{\mathcal A}}
\newcommand{\calB}{{\mathcal B}}
\newcommand{\calC}{{\mathcal C}}
\newcommand{\calD}{{\mathcal D}}
\newcommand{\calF}{{\mathcal F}}
\newcommand{\calI}{{\mathcal I}}
\newcommand{\calK}{{\mathcal K}}
\newcommand{\calM}{{\mathcal M}}
\newcommand{\calP}{{\mathcal P}}
\newcommand{\calR}{{\mathcal R}} 
\newcommand{\calS}{{\mathcal S}}
\newcommand{\calT}{{\mathcal T}}
\newcommand{\calU}{{\mathcal U}}
\newcommand{\calV}{{\mathcal V}}
\newcommand{\calW}{{\mathcal W}}
\newcommand{\spec}{\mbox{spec}\,}
\newcommand{\phg}{{\mathrm{phg}}\,}
\newcommand{\coker}{\mathrm{coker}\,}
\newcommand{\MM}{{\mathbb M}}

\newtheorem{theorem}{Theorem}
\newtheorem{proposition}{Proposition}
\newtheorem{corollary}{Corollary}
\newtheorem{lemma}{Lemma}
\newtheorem{definition}{Definition}
\newtheorem{remark}{Remark}

\begin{abstract} 
The theory of complete surfaces of (nonzero) constant mean curvature
in $\RR^3$ has progressed markedly in the last decade. This paper surveys
a number of these developments in the setting of Alexandrov embedded surfaces;
the focus is on gluing constructions and moduli space theory, and the analytic 
techniques on which these results depend. The last section contains some
new results about smoothing the moduli space and about CMC surfaces in
asymptotically Euclidean manifolds.
\end{abstract}

\section{Introduction}
The study of global minimal and constant mean curvature (CMC) surfaces in $\RR^3$ has always 
occupied a central position in classical differential geometry, but the former class of surfaces
has undoubtedly received more attention than the latter. While partly historical accident, 
this is also due to the classical Weierstra{\ss} representation formula for minimal
surfaces, through which many interesting examples of minimal surfaces were known.
In contrast, the most famous results about CMC surfaces prior to 1980 were the rigidity theorems 
of Hopf and Alexandrov, which state that the only compact CMC surfaces in $\RR^3$ which
are immersed and genus $0$, respectively embedded of arbitrary genus, are spheres; 
in some sense this is a negative result. There has been an explosive development in both these fields
during the past twenty years, centering on many new examples and methods of construction,
and also including various classification and structure theorems. The recent volume \cite{Clay} 
contains up-to-the-minute accounts of many of these developments, but mostly focussing on minimal
surfaces. My goal here is to survey some of the advances in the global study of CMC surfaces
from my own biased perspective. 

Let $\Si$ be a surface in $\RR^3$. Its mean curvature $H$ at a point $p$ is the sum of the two 
principal curvatures (the customary factor of $1/2$ is omitted), and $\Si$ has constant 
mean curvature if this quantity is constant. In this case we can normalize so that 
either $H \equiv 0$ or $H \equiv 1$. We shall restrict attention exclusively to the second case.
In addition, we shall only be concerned with complete, properly immersed surfaces of finite
topology. There is a rich theory for the Plateau problem for minimal and CMC surfaces, 
cf.\ \cite{Struwe}, \cite{BC}, but we shall not discuss this at all; in particular, we
shall understand all surfaces below to be complete, unless otherwise stated. On the other hand,
unlike the situation for minimal surfaces, cf. \cite{CKMR} for example, the subject of CMC 
surfaces with infinite topology is mostly unexplored.

The only complete CMC surfaces known classically are the sphere, the cylinder, and 
the one-parameter family of rotationally invariant Delaunay surfaces, which we describe in the 
next section. (We should also mention a method, due originally to Bonnet but most clearly
enunciated by Lawson, for associating a `cousin' CMC surface in $\RR^3$ to a minimal surface in 
$\sphere^3$.) The first modern breakthrough was Wente's discovery, detailed in his 1983 paper \cite{W},
of immersed CMC tori. Not long afterwards, Kapouleas used transcendental PDE 
methods to construct many new CMC surfaces, including compact ones with arbitrary genus
\cite{Ka2}, \cite{Ka3}, and noncompact ones with finitely many ends \cite{Ka1}. 
Gro{\ss}e-Brauckman \cite{KGB} subsequently constructed certain of these noncompact 
surfaces with large discrete symmetry groups using more classical methods based on Schwarz reflection.

The theory has developed in two fairly distinct directions. From the techniques used for
Wente's construction ultimately has emerged the DPW (Dorfmeister-Pedit-Wu) method, which 
draws on the theory of integrable systems, and serves as a replacement for the
Weierstra{\ss} representation; on the other hand, Kapouleas' construction has engendered 
many new PDE approaches to the problem. These 
theories have distinct flavors, and in some senses illuminate quite different 
classes of surfaces. For example, the former method seems to be more successful in describing
general immersed CMC surfaces, and is closely related to many of the computer experiments
and simulations of CMC surfaces, cf.\ \cite{KMS}, while the latter has been particularly good 
in describing CMC surfaces which are `nearly embedded'. We specialize once again, and for the 
last time, to the second PDE set of methods. 

It turns out that requiring CMC surfaces to be embedded is too restrictive, 
and in some sense even geometrically unnatural. Thus, we shall relax this hypothesis and 
require instead that our CMC surfaces be {\it Alexandrov embedded}. To define this, first 
recall that a complete CMC surface $\Si$ with 
finite topology is diffeomorphic to a punctured Riemann surface $\overline{\Si}\setminus
\{p_1, \ldots, p_k\}$; we shall always assume that $\Si$ and $\overline{\Si}$ are oriented. 
Writing $\overline{\Si} = \del Y^3$, then we say that $\Si$ is Alexandrov embedded
if the immersion $\Si \hookrightarrow \RR^3$ extends to an immersion $Y \hookrightarrow \RR^3$. 
This condition is naturally suited to the technique of Alexandrov reflection, cf.\ \cite{KKS}. 

The basic questions we shall address are whether there are any Alexandrov embedded CMC surfaces,
and if so, how might we go about constructing them and describing the totality of them. 
The next section discusses the main elementary examples and the basic structure theory of 
complete Alexandrov embedded CMC surfaces, and lays out the framework for the subsequent 
develoment. In \S 3 we describe a variety of analytic constructions, centered on the
basic technique of Cauchy data matching. This is followed, in \S 4, by the ramifications
of these constructions for the moduli space theory. Finally, in \S 5, we generalize 
this theory to CMC surfaces in arbitrary asymptotically Euclidean $3$-manifolds and
suggest some possible directions of future research.

This survey has some overlap with one given in the recent paper by Kusner \cite{Ku}, 
which I recommend strongly to the interested reader; his emphasis is rather different 
than mine and I hope that these two papers complement one another in a useful way.
I should also mention that the forthcoming thesis of M.\ Jleli at Universit\'e Paris XII 
contains the analogues of essentially all of the results here for CMC hypersurfaces in $\RR^n$. 

My papers in this area have been written in various collaborations 
with Rob Kusner, Frank Pacard, Dan Pollack and Jesse Ratzkin. I owe special thanks to Frank 
Pacard for all he has taught me over the years! Kusner, Pacard and Pollack each
read and gave me many helpful comments on a preliminary draft of this paper, 
but I happily take credit for the 
remaining errors and solecisms. As described in a previous survey \cite{MPo}, most of the 
results in this theory have parallels for solutions of the singular Yamabe problem with isolated 
singularities. In particular, the general moduli space results in \S 4 originated 
in my work with Pollack and Uhlenbeck \cite{MPU} on deformations of solutions of the Yamabe
equation with isolated singularities on the sphere. I am also grateful to Nick Korevaar
and Rick Schoen for many illuminating discussions. I thank the Mittag-Leffler
Institute and Universita di Roma ``La Sapienza'' for their hospitality while this paper was written.
Finally, this paper is based on a talk given at the conference in honour of H.\ Brezis and
F.\ Browder at Rutgers University in October 2001, and I would like to acknowledge that 
my first exposure to CMC surfaces was in a course taught by Professor Brezis at MIT in 
the early 1980's. 

\section{Basic properties}
In this section we introduce most of the basic definitions. Here and in the
remainder of this paper we denote by $\calM_{g,k}$ the set of all
proper Alexandrov embedded constant mean curvature $1$ surfaces in $\RR^3$ of genus $g$ with $k$ ends, 
and when we say that a surface $\Si$ is CMC, we mean that it lies in one of these
spaces. We begin by discussing the most elementary examples of complete noncompact 
CMC surfaces, the Delaunay surfaces. This is followed by a description of the 
basic structure theory of (Alexandrov embedded) CMC surfaces. We go on to discuss the Jacobi 
operator and its mapping properties, and the fundamental notion of nondegeneracy of a CMC surface. 

\subsection{Delaunay surfaces} The CMC surface with the largest symmetry group is, of course, 
the sphere, but the most familiar one after that is the cylinder. The sphere is the only 
surface (CMC or not) with full rotational symmetry, and by the classical results of Hopf and 
Alexandrov mentioned earlier, it is the only compact Alexandrov embedded CMC surface; on the other 
hand, there is an interesting
family of axially symmetric CMC surfaces discovered by Delaunay \cite{D}, in which
the cylinder lies as an extreme element. To describe these, consider the cylindrical
graph $(t,\theta) \mapsto (\rho(t)\cos\theta, \rho(t)\sin\theta,t)$. The CMC equation
for this graph is an ODE, 
\begin{equation}
\rho_{tt} - \frac{1}{\rho}(1+\rho_t^2) + (1 + \rho_t^2)^{3/2} = 0.
\label{eq:deleqn}
\end{equation}
All positive solutions of this equation are periodic; this follows from the facts that the function 
$H(\rho,\rho_t) = \rho^2 - 2\rho(1+\rho_t^2)^{-1/2}$ is an integral of motion of this ODE and 
has compact level curves in $\{\rho>0, H < 0\}$. (Kusner's balancing formula, described below, gives
a nice geometric interpretation of the fact that $H$ is constant on solutions.) 
For any $0 < \e < 1$ let $\rho_\e(t)$ be the solution which attains the minimum value $\e$, and 
normalize by translating the independent variable so that $\rho_\e(0)=\e$. 
Then $\e \leq \rho_\e(t) \leq 2-\e$ for all $t$, and 
since $\rho_\e$ remains positive, the surfaces $D_\e$ which are the cylindrical graphs of 
these functions are embedded; they are called Delaunay unduloids. 

There are two limiting cases: in the first, $\e=1$ and $D_1$ is the cylinder of radius $1$; 
the other occurs when $\e \searrow 0$, and then $D_\e$ converges to an infinite array 
of mutually tangent spheres of radius $2$ with centers along the $z$-axis. The family $D_\e$ 
interpolates between these two extremes. The number $\e$ measures the size of the `neck' region. 

The Delaunay family continues past the singular limit at $\e=0$ to a family of immersed 
(not Alexandrov embedded!) CMC surfaces known as the Delaunay nodoids. These will only
play a minor role in this survey. 

There is an alternate, conformal, parametrization of the $D_\e$ which is more convenient for 
most calculations. This isothermal parametrization involves two functions $\sigma$ and $\kappa$:  
$\sigma$ is the unique smooth nonconstant solution to the initial value problem
\[
\left(\frac{d\sigma}{ds}\right)^2 + \tau^2 \cosh^2 \sigma =1, \qquad 
\del_s\sigma(0) = 0, \quad \sigma (0) < 0, \quad \mbox{when}\ 
\tau \in (0,1]
\]
while $\kappa(s)$ is determined by
\[
\frac{d\kappa}{ds} = \ \tau^2 \, e^{\sigma}\, \cosh \sigma, \qquad 
\kappa(0) = 0.
\]
Setting $t = \kappa(s)$ and $\rho (\kappa(s)) = \tau e^{\sigma(s)}$, 
and determining $\tau$ by $\e = 1 - \sqrt{1-\tau^2}$, then $D_\e$ is parametrized by 
\begin{equation}
X_\e: \RR \times S^1 \ni  (s,\theta) \mapsto 
\left(\tau\,e^{\sigma_\tau(s)} \, \cos\theta,\tau\,
e^{\sigma_\tau(s)}\,\sin\theta,\kappa_\tau (s)\right).
\label{eq:1.1}
\end{equation}

The metric coefficients in this coordinate system are $g_{ss} = 
g_{\theta \theta} = \tau^2 e^{2\sigma}$, and $g_{s \theta}= g_{\theta s} = 0$. 

While this definition may seem ad hoc, it is motivated by a systematic line of 
reasoning in the theory of integrable systems; \cite{MP1} contains a more
complete discussion of this parametrization. 

\subsection{Structure theory} There are several important 
facts about the global structure of Alexandrov-embedded CMC surfaces
of finite topology, and we review these now.

The story begins with an important theorem of Meeks \cite{Me}, which states that 
any proper Alexandrov embedded end of a CMC surface is cylindrically bounded, i.e.\ 
contained in some tube of radius $R>0$ around a ray. This result reflects
the `smallness at infinity' of the ambient Euclidean space and fails definitively,
for example, in hyperbolic space. (Simple counterexamples are the equidistant
surfaces around a totally geodesic copy of $\HH^2 \subset \HH^3$.)  Meeks also
proves that there are no one-ended CMC surfaces, i.e.\ $\calM_{g,1} = \emptyset$.

Building on this, Korevaar, Kusner and Solomon \cite{KKS} obtained a much
sharper picture of the behaviour of CMC ends. More precisely, for any such end $E$ there 
is a half Delaunay unduloid $D_\e^+ = D_\e \cap \{z \geq 0\}$ and a rigid
motion $A$ of $\RR^3$ such that $E$ converges exponentially to $A(D_\e^+)$.
More precisely, if $\nu$ is the outward unit normal to $A(D_\e^+)$, then any
surface $\calC^1$ close to this Delaunay end can be written as a normal graph, 
i.e.\ as the image of 
\[
A(D_\e^+) \ni p \longrightarrow  p + \phi(p)\nu(p).
\]
Then $E$ is the normal graph of a function $\phi \in \calC^\infty(A(D_\e^+))$
which decays exponentially at infinity. The proof ultimately devolves to
the exponential decay of Jacobi fields (see below) on a Delaunay surface
and an adaptation of a beautiful `approximation improvement' lemma due to L.\ Simon.
Using Alexandrov reflection they prove that $\calM_{g,2} = \emptyset$ if $g \geq 1$, while 
$\calM_{0,2}$ contains only the Delaunay surfaces. 

Amongst the many corollaries of this result is the fact that we can associate
an asymptotic  necksize parameter $\e$ to any end of $\Si \in \calM_{g,k}$. 
There are other ways to measure this necksize. For example, there is a 
`balancing formula' for CMC surfaces, discovered by Kusner,  
which states the following. Let $\Gamma$ be a simple closed curve on $\Si$ and $C$ an
immersed surface in $\RR^3$ with $\del C = \Gamma$; suppose that $n$ and $\nu$ are the unit normals 
to $\Gamma$ in $T\Si$ and $C$, respectively. (Thus $n$ is orthogonal to the unit normal 
to $\Si$ in $\RR^3$.) Then, with an appropriate choice of orientations of these normals, 
the quantity
\[
\frac{1}{\pi}\left(\int_\Gamma n \, d\sigma - \int_C \nu \, dA \right)
\]
depends only on the homology class of $\Gamma$ in $\Si$. Thus, for example, if $\Gamma$ 
lies in an end $E$ of $\Si$ and represents a generator of $H_1(E)$, then it is
homologous to the curve $\Gamma_s$ which is the normal graph in $E$ of the circle
$\theta \to X_\e(s,\theta)$ in the model Delaunay end, for any $s \geq s_0$.
Using the exponential decay of $E$ to this Delaunay end, this quantity can be 
computed explicitly, and in fact just equals $\e(2 - \e)\vec v$,
where $\vec v$ is the unit vector in the direction of the Delaunay axis. 
This is called a balancing formula because if $\Gamma_1, \ldots, \Gamma_k$
encircle the $k$ ends of $\Si \in \calM_{g,k}$, and are oriented so that their
sum is null homologous, then the sum of the corresponding force vectors must vanish. 

In a subsequent paper \cite{KKMS}, these authors and Meeks establish an analog of this 
asymptotics result for proper Alexandrov embedded CMC surfaces in the hyperbolic space  
$\HH^3$, in the case that the mean curvature $H>2$. We note that CMC (hypersurface)
theory in $\HH^n$ is quite different in the three cases, $H>n-1$, $H<n-1$ and $H=n-1$, 
the latter having the most similarities with minimal surface theory (via the cousin 
correspondence mentioned earlier), and the first having the most similarities with
CMC theory in $\RR^n$. \cite{KKMS} also gives Delauany asymptotics for ends of CMC 
hypersurfaces in $\HH^n$, but only under the {\it a priori} assumption of cylindrical 
boundedness. There is a growing literature on CMC surfaces in hyperbolic space, 
cf.\ \cite{Br}, \cite{Pa}, \cite{GM}, and the 
work of H.\ Rosenberg and his collaborators.
It is quite likely that the results of \cite{KKS} and \cite{KKMS}
remain true for CMC surfaces in manifolds which are only asymptotically Euclidean
or hyperbolic, and indeed, the proofs from these papers may require only
moderate alterations.  We come back to this point in \S 5. 

These asymptotics results definitely require Alexandrov-embeddness,
or something near to it. For example, the Delaunay family $D_\e$ continues
beyond $\e=0$ to a family of immersed (not Alexandrov embedded!) rotationally
invariant CMC surfaces called nodoids, and Pacard and I show \cite{MP2} that 
these become increasingly unstable as $\e \to -\infty$: there are infinitely
many values $\e_j \to -\infty$ where new CMC surfaces which are not rotationally
symmetric bifurcate away from this nodoid family. This bifurcation bears some relationship
to the classical Rayleigh instability of the cylinder. It only
occurs when $\e \leq \e_0 < 0$, and it is possible that there may be some
global structure theory, akin to that described above, for properly immersed CMC 
surfaces for which all ends have force vector corresponding to some Delaunay
nodoid $D_\e$ with $\e_0 < \e < 0$. We leave this as an interesting 
unexplored direction. 

One can interpret this asymptotics theorem by thinking of $\Si \in \calM_{g,k}$ as decomposing 
into a large compact piece, about which we know very little {\it a priori}, with $k$ 
asymptotically Delaunay ends emerging from it. All of the gluing constructions described
in \S 3 involve ``filling in the black box in the middle'' in different ways. Kapouleas'
original construction builds CMC surfaces around (suitably balanced) simplicial graphs where 
$k$ of the edges are rays tending to infinity. The finite edges and the rays are replaced 
by finite and semi-infinite segments of Delaunay surfaces, respectively, and the vertices
become spheres. (Balancing here means that there are force vectors associated
with each edge and these must cancel at each vertex. An additional condition regarding
the flexibility of this arrangement is also needed.) 
As we discuss later, there are various other geometric possibilities
for the vertices, but Korevaar and Kusner \cite{KK} show that this picture is
qualitatively correct in general. For any $\Si \in \calM_{g,k}$, there is a
simplicial graph, with explicit bounds on the numbers of edges and vertices in terms 
of $g$ and $k$, and a (fairly large) tubular neighbourhood around it which
contains $\Si$. 

\subsection{The Jacobi operator: its mapping properties, nullspace and nondegeneracy}
If $\Si$ is a complete Alexandrov embedded CMC surface with finite topology, 
then we can consider all surfaces which are $\calC^2$ close to $\Si$ as normal
graphs, i.e.\ for any $\phi \in \calC^2(\Si)$, sufficiently small, we set
\begin{equation}
\Si_\phi = \{x + \phi(x)\nu(x): x \in \Si\}
\label{eq:normalgraph}
\end{equation}
where $\nu(x)$ is the unit normal to $\Si$ at $x$. The equation which determines that
$\Si_\phi$ also has CMC equal to $1$ is a quasilinear elliptic equation for $\phi$.
Its linearization at $\phi = 0$ is called the Jacobi operator $L$ (or $L_\Si$). It
is well known that
\[
L_\Si = \Delta_\Si + |A_\Si|^2,
\]
where $A_\Si$ is the second fundamental form. Solutions of $L \phi = 0$ are
called Jacobi fields. If $u_t$ is a $1$-parameter family of functions such
that all of the surfaces $\Si_{u_t}$ are all CMC (with the same mean curvature),
then $\phi = \left. du_t/dt\right|_{t=0}$ is a Jacobi field on $\Si$. In
general, Jacobi fields need only correspond to deformations of $\Si$ which are 
CMC to second order. 

One very important consequence of the structure theory described in the last 
subsection is that because CMC surfaces have very well-controlled geometry
at infinity, the behaviour of their Jacobi operators is governed by the
behaviour of Jacobi operators on (half) Delaunay surfaces. Indeed, as we
discuss in a moment, the asymptotics of any (tempered) Jacobi field on $\Si$
are determined by the asymptotics of tempered Jacobi fields on the model
Delaunay ends. More generally, mapping properties for $L_\Si$
may be deduced from those for the Jacobi operators $L_\e$ on $D_\e$. 
This motivates a closer study of these model operators $L_\e$. 

In cylindrical coordinates the expression for $L_\e$ is quite complicated, but 
using the isothermal parametrization, the Jacobi operator $L_\e$ on $D_\e$ takes the 
much simpler form
\[
L_\e = \frac{1}{\tau^2e^{2\sigma}}\left(\del_{s}^2 +  \del_{\theta}^2 + 
\tau^2 \, \cosh (2\sigma) \right). 
\]
Since $\tau^2 e^{2\sigma}$ is bounded away from zero and infinity, it is
sufficient to analyze the operator
\[
\del^2_{s} + \del^2_{\theta} + \tau^2 \, \cosh (2\sigma(s)), 
\]
and this is what we do in practice. This analysis relies on the 
rotational invariance of this operator in $\theta$ and its periodicity 
in $s$. Thus we can simultaneously reduce by Fourier
series in the $\sphere^1$ factor and use an adapted form of Flocquet
theory to handle the resulting ODEs with periodic coefficients. 
With these techniques we are able to determine that $L_\e$ is Fredholm
on certain natural exponentially weighted function spaces, and to characterize
all global tempered solutions. We return to this last point shortly.
The analysis of $L_\Si$ itself is accomplished by patching the parametrices
for each of these model problems, corresponding to each of the ends $E_j$
of $\Si$, together with an interior parametrix. These arguments are presented 
in detail (for a slightly different geometric problem) in \cite{MPU}, cf.\ also 
\cite{MPP2}. 

We shall discuss Jacobi fields on $\Si$ in more detail in \S 4. For the present
we note that one may generate a distinguished set of Jacobi fields on any CMC 
surface using rigid motions of the ambient space. The corresponding 
`geometric Jacobi fields' thus correspond to infinitesimal translations and rotations 
of $\Si$ in $\RR^3$. For a Delaunay surface $D_\e$, the three independent
translations induce three Jacobi fields, and the two infinitesimal rotations of the
axis of symmetry induce two more. $D_\e$ has a final Jacobi field, induced 
by change of the necksize parameter $\e$. (The case $\e=1$ is handled slightly
differently.) These all have the form $e^{i\ell\theta}\psi(s)$, where $\ell = 0, \pm 1$ (or rather, 
are real-valued combinations of these); for the Jacobi fields coming from 
translations, $\psi(s)$ is periodic, while for each of the others it grows linearly as 
$s \to \pm \infty$. These six Jacobi fields are the only global solutions of $L_\e \phi = 0$ on
$D_\e$ which do not grow exponentially on one end or the other. In fact, there
are solutions of the form $e^{i\ell\theta}\psi(s)$ for any $\ell \in \ZZ$, $|\ell| \geq 2$,
but for each of these, $\psi(s+S_\e) = e^{\gamma_\ell}\psi(s)$, where $S_\e$ is the period
of $\sigma$ and where $\gamma_\ell \in \RR$, $\gamma_\ell \neq 0$. In particular, these
other solutions have precise exponential rates of growth or decay. 
This set of exponents constitutes a discrete set $\Lambda_\e$, called the
indicial set of $L_\e$. 

As for the mapping properties of $L_\Si$, one's first temptation might be to
let this operator act on $L^2(\Si)$. However, this mapping does not have closed
range. In fact, the parametrix method sketched above shows that the spectrum
of $L_\Si$ on $L^2$ is comprised of bands of continuous spectrum and locally 
finite discrete spectrum. Unfortunately, $0$ is {\it always} on the end of
one of these bands of continuous spectrum. In the language of scattering theory,
we say that it is a threshold value. To see this, recall from above that $\Si$ always 
admits nontrivial global bounded and linearly growing Jacobi fields on $\Si$, 
corresponding to translations and rotations in $\RR^3$.  This is already enough
to imply that $0$ is in the continuous spectrum by Weyl's criterion; 
a closer examination of the Flocquet theory, together with the existence of 
linearly growing (as opposed to bounded, periodic) Jacobi fields, shows that $0$ 
is a threshold value. 

In any case, one must look elsewhere to realize $L_\Si$ as an
operator with closed range. As indicated earlier, this may be done using 
spaces with exponential weights. We can use either Sobolev or H\"older spaces,
and each has its advantages, but we shall use the latter. 
Decompose $\Si$ into the union of a compact
piece, $K$, and a finite number of ends, $E_1, \ldots, E_k$. Fix isothermal
coordinates $(s_j,\theta_j)$ on each end. 
\begin{definition} For $\ell \in \NN$, $0 < \al < 1$ and $\mu \in \RR$, we define 
$\calC^{\ell,\al}_\mu(\Si)$ to be the set of all functions $u\in \calC^{\ell,\al}(K)$ and
on each end $E_j$ satisfy 
\[
\sup_{s_j \geq 0} e^{-\mu s_j}\, \| w\|_{\calC^{\ell,\alpha}([s_j,s_j + 1]\times S^1)} < \infty.
\]
\end{definition}
It is immediate that 
\begin{equation}
L: \calC^{\ell+2,\al}_\mu(\Si) \longrightarrow \calC^{\ell,\al}_\mu(\Si)
\label{eq:mL}
\end{equation}
is a bounded mapping for any $\mu \in \RR$. However, if $\e_j$ is the
asymptotic necksize of the end $E_j$, then it is not hard to see that
whenever $\mu \in \Lambda_{\e_j}$, (\ref{eq:mL}) does not have closed range.
However, this is the only bad case:
\begin{proposition} Let $\Lambda = \Lambda_\Si = \cup_{j=1}^k \Lambda_{\e_j}$. Then if $\mu \notin 
\Lambda$, the mapping (\ref{eq:mL}) is Fredholm. When $\mu \notin \Lambda$ is sufficiently large, 
this mapping is surjective, while if $\mu$ is sufficiently negative, this mapping is injective.
\end{proposition}

We conclude this subsection by stating the fundamental
\begin{definition} The CMC surface $\Si$ is said to be {\bf nondegenerate}
if it admits no nontrivial global Jacobi fields which decay along all of the
ends $E_j$, or equivalently, if (\ref{eq:mL}) is injective for any $\mu < 0$.
\end{definition}

The equivalence of the two versions of this definition is not immediate,
since one could imagine the existence of Jacobi fields which decay polynomially 
rather than exponentially, 
but is true nonetheless.  A duality argument then gives that if $\Si$
is nondegenerate, (\ref{eq:mL}) is surjective when $\mu > 0$, $\mu \notin \Lambda$.
There is an important refinement of this surjectivity. Choose a cutoff
function $\chi(s)$ which equals $1$ for $s \geq 1$ and $0$ for $s \leq 0$,
and consider the $6k$-dimensional vector space $W_\Si$ spanned by the elements
\[
\Phi_{ij} = \chi(s_j)\phi_{i,\e_j}(s_j,\theta_j),\ j = 1, \ldots, k,\ i = 1, \ldots, 6,
\]
where $\{\phi_{i,\e_j}\}$ are the $6$ geometric Jacobi fields on $D_{\e_j}$. Thus we have 
transplanted these special Jacobi fields to each end of $\Si$. Because of the exponential
decay of the ends to their Delaunay models, we have
\[
L_\Si \, \phi = {\mathcal O}(e^{-c s_j}) \qquad \forall \phi \in W_\Si,\ j = 1, \ldots, k.
\]
Here $c = \inf\{\gamma: 0 < \gamma \in \Lambda_\Si\} > 0$.  
Hence for $\mu > 0$ sufficiently small,
\begin{equation}
L_\Si: W_\Si \oplus \calC^{k+2,\al}_{-\mu}(\Si) \longrightarrow \calC^{k,\al}_{-\mu}(\Si) 
\label{eq:aug}
\end{equation}
is well defined.
\begin{proposition} If $\Si$ is nondegenerate and $0 < \mu < 
\mu_0 = \inf \{|\lambda|: \lambda \in \Lambda_{\Si}\}$, 
then the mapping (\ref{eq:aug}) is surjective.
\end{proposition}

We call $W_\Si$ the deficiency subspace of $\Si$. It also plays a role in the 
regularity of Jacobi fields. 
\begin{proposition} Let $\phi$ be any Jacobi field with at most polynomial growth on $\Si$.
Then there exist functions $\Phi \in W_\Si$ and $\psi \in \calC^{2,\al}_{-\mu}(\Si)$,
$\mu > 0$, such that $\phi = \Phi + \psi$.
\end{proposition}

The nonlinear mean curvature operator $N$ does not preserve any of the spaces $\calC^{k,\al}_\mu(\Si)$, 
$\mu > 0$, because the nonlinearity amplifies the exponential increase.
On the other hand, while this nonlinear operator is defined on the spaces of
exponentially decreasing functions $\calC^{k,\al}_{-\mu}$, its linearization always has a 
cokernel there. We can extend the definition of $N$ to $W_\Si$ as follows: any
$\Phi \in W_\Si$ corresponds to a one-parameter family of CMC deformations of the
model Delaunay surface for that end, and so we can `integrate' $\Phi$ by transplanting this
to a one-parameter family of asymptotically CMC deformations localized to each of
the ends of $\Si$. For example, if $\Phi$ corresponds to the rotation of a Delaunay
surface, then we rotate the corresponding end, using some cutoff to join this to
the identity on the rest of $\Si$. The constant mean curvature condition is destroyed
only at the pivoting locus which is contained in a compact set. (The definition
is slightly more complicated when $\Phi$ includes a component requiring a change
of Delaunay parameter.) More generally, $N(\Phi)$ 
is just the mean curvature of this new surface, and this is either compactly
supported or exponentially decreasing, so that
\begin{equation}
N: W_\Si \oplus \calC^{k+2,\al}_{-\mu}(\Si) \longrightarrow \calC^{k,\al}_{-\mu}(\Si) 
\label{eq:augnl}
\end{equation}
is well defined. When $\Si$ is nondegenerate, this map has surjective
linearization. Thus the addition of $W_\Si$ provides an intermediate space where the 
conflicting requirements of well-definedness and surjectivity are balanced.  

We conclude this section by noting that there is a more general way to define 
nondegeneracy. The crucial feature in the discussion above is that we were
able to extend the definition of $N$ to $W_\Si$ because elements of the deficiency
space correspond, at least asymptotically, to actual CMC deformations. Another
way to say this is that every element of $W_\Si$ asymptotically integrates to
a one-parameter family of asymptotically CMC surfaces. Thus the `correct' 
definition of nondegeneracy should be to require that all tempered 
Jacobi fields on $\Si$ are similarly asymptotically integrable. The point
here is that one should be able to use explicit geometric deformations of $\Si$
to compensate for the cokernel of $L_\Si$. As a very
simple but important example, the sphere $\sphere^2$ is degenerate in the
more restricted sense; it has a three-dimensional space of Jacobi fields,
which are the restrictions of the linear coordinates on $\RR^3$ to $\sphere^2$. 
However, these Jacobi fields correspond to translations of the sphere, and
it is precisely this which allowed Kapouleas to use them in his gluing
construction. 

The possible existence of degenerate CMC surfaces creates a lot of complications
throughout this whole theory. However, even amongst immersed surfaces, there
are no known examples which are degenerate in the sense of this broader
definition. Indeed the only known examples of degenerate surfaces in the narrower
sense are the sphere (or any compact immersed CMC surface) and certain immersed 
`bubbleton' solutions which are globally cylindrically bounded and asymptotically 
cylindrical, so that translation along the cylinder's axis provides a decaying, but nevertheless
`integrable', Jacobi field. A major problem is to determine whether there
ever exists degenerate CMC surfaces in the broader sense. 

\section{Gluing constructions} Up to this point the only examples of complete Alexandrov 
embedded CMC surfaces of finite topology we have seen are the Delaunay unduloids. In this 
section we discuss a number of closely related gluing constructions; taken together 
these show that CMC surfaces are quite malleable and exist in great profusion. 
As with all gluing constructions, the rough idea is that one starts with 
`building blocks', pieces of surfaces which already are, or are near to, 
simpler CMC surfaces, and then pieces them together to form complete CMC surfaces. 
We use here three types of building blocks: Delaunay surfaces (or finite
or semi-infinite segments of them), spheres and complete Alexandrov-embedded 
{\it minimal} surfaces with finite total curvature. It is necessary that
all components be nondegenerate in the extended sense. Amongst the motivations 
for developing new analytic methods for these sorts of problems
is that one wishes to show that the glued surfaces are again nondegenerate;
this is important in the moduli space theory and also has the important
consequence that these constructions can be iterated. 
In addition, until the advent of newer methods developed 
after Kapouleas' work, it was unclear whether there existed any nondegenerate 
surfaces beyond the unduloids. 

I shall first give a rough overview of the type of argument needed here
(or indeed any singular perturbation problem). After describing Kapouleas' 
construction I list the more general gluing constructions now known to
work, and conclude by describing idea of Cauchy data matching which can
be used to prove any of these.  All of these surfaces have very small necksizes,
or otherwise nearly degenerate (i.e.\ near to the boundary of Teichm\"uller space)
conformal structures. It is an important philosophical
point that all of the `constructible' CMC surfaces are very nearly degenerate.
Gluing constructions provide parametrizations of ends of the moduli spaces.
The (presumably more commonplace) surfaces with larger necksizes are not
`reachable' by direct methods. 

Thus, again at a rough level, these components are arranged, by rotating and translating 
them in space, so that they form a configuration which is an approximate solution. 
Such a configuration is either CMC or approximately CMC everywhere. There 
is always a parameter $\eta$ involved, so that as $\eta \searrow 0$, the discrepancy of 
the configuration from being exactly CMC diminishes. Unfortunately, in this same limit, 
the geometry of the configuration degenerates. The main step involves using the 
Jacobi operator to perturb the approximate configuration to an exact CMC surface,
but again this is an operator with `degenerating coefficients'.  Hence we must
not only show that its Jacobi operator is surjective as in (\ref{eq:aug}) when $\eta$ 
is small, but also estimate the growth of the norm of a suitable right inverse as
$\eta \searrow 0$. The game is to show that this norm does not blow up faster than 
the inverse of the size of the error term.  There are many well-known constructions
which follow these outlines, including Taubes' famous instanton patching, etc., 
and there are many different methods of carrying out the details. 
I shall sketch one, initially developed in \cite{MP1}, based on Cauchy data 
matching; this seems to be the most efficient and involves the fewest hard estimates. 

The first construction of this type for CMC surfaces was accomplished by 
Kapouleas \cite{Ka1}, inspired by a similar construction by Schoen \cite{Sch} for the 
singular Yamabe problem. This 
establishes the existence of many complete CMC surfaces, but his method does not give
a sufficiently fine understanding of their geometry; in particular, it seems very 
difficult to tell whether any of the surfaces he constructs are nondegenerate.
(They are, in fact, nondegenerate, as follows from the alternate 
construction v) below.) 
I have already sketched his method in \S 2: one begins with a simplicial graph 
in $\RR^3$ with $k$ of the edges semi-infinite rays, and where each edge is labelled
with a force vector. This graph must satisfy a certain balancing and flexibility
condition. The approximate CMC surface is assembled by substituting for each
vertex a sphere, for each finite edge a finite Delaunay segment (possibly just a single
neck region), and for each ray a half Delaunay surface. The force vectors include
information about the relative scaling of the Delaunay necksizes, and the
parameter $\eta$ gives the absolute size of these necks. Kapouleas joins 
these pieces with cutoff functions; the resulting error terms are large in $L^\infty$, 
but are supported in small sets. A delicate analysis of the Jacobi operator 
(complicated by the fact that he works on a space where the Laplacian does not have
closed range), allows him to show that the necessary perturbation can be made 
when $\eta$ is small.

We now list and give brief descriptions of various geometric operations 
which, it has more recently been proved, may all be done in the (nondegenerate)
CMC category. Afterwards we sketch a representative proof.

\begin{itemize}
\item[i)] Adding Delaunay ends \cite{MPP2}: If $\Si$ is any nondegenerate CMC surface and
$p \in \Si$, then we may attach a half Delaunay surface $D_\e^+$ with 
small necksize $\e$ to $\Si$ at $p$. The axis of $D_\e^+$ is directed in
the outward normal direction. 
\item[ii)] Connected sums \cite{MPP1}: Let $\Si_1$ and $\Si_2$ be nondegenerate CMC
surfaces and choose points $p_j \in \Si_j$. Rotate and translate these
surfaces so that the tangent spaces coincide at these points, but with
opposite orientation. Then there exists a new CMC surface $\Si_1 \#_\e \Si_2$
obtained by `bridging' these surfaces with a small Delaunay neck near
these points. (One can also insert any finite Delaunay segment instead of a single neck.) 
To show that the resulting surface is nondegenerate requires
minor restrictions concerning the locations of the points.
\item[iii)] Attaching Delaunay ends to minimal $k$-noids \cite{MP1}: Let $\Si$
be a complete nondegenerate {\it minimal} surface of finite total curvature with $k$ 
ends (a so-called minimal $k$-noid), all of which are asymptotic to catenoids. Form a 
new surface with boundary $\Si_\e$ by truncating the ends (at distance proportional to 
$-1/\log \e$ on the catenoid) and then rescaling by a factor of $\e$. Then there is a 
CMC surface with $k$ ends obtained by attaching half-Delaunay surfaces 
$D_{\e_j}^+$ to these truncated catenoidal ends. The necksize parameters are 
$\e a_j + {\mathcal O}(\e^2)$, where $(a_1, \ldots, a_k)$ is a vector of relative dilation 
factors for the catenoids modelling the ends of $\Si$. Since many topologically
complicated nondegenerate minimal $k$-noids are known to exist, this gives similarly
complex CMC surfaces. It may seem that being {\it minimal}
is very far away from having constant mean curvature $1$, but this is where the 
scaling by $\e$ becomes important. For, rescaling the ensemble by $1/\e$, the construction
is equivalent to attaching (albeit on a very large scale) Delaunay ends with mean 
curvature $\e$ to a minimal surface, which seems more plausible.
\item[iv)] End-to-end gluing \cite{Rat}: Let $\Si_1$ and $\Si_2$ be two
nondegenerate CMC surfaces, each of which has an end with the same asymptotic
necksize parameter (not necessarily small). Then we may translate and rotate
these surfaces so that the axes of these particular ends are on the same
line, but oppositely oriented. Truncating these ends very far out, it is 
possible to attach the surfaces to obtain a new CMC surface with a very 
long approximately unduloidal tube. This requires a small strengthening
of the nondegeneracy condition. Notice that while these surfaces do not necessarily
have small necksizes, they contain embedded essential annuli with large conformal
modulus, hence still have nearly degenerate conformal structures.
\item[v)] Kapouleas-type constructions \cite{MPPR}: One may attach half Delaunay 
surfaces and finite Delaunay segments, arrayed along the rays and edges of a
suitable simplicial graph, with spheres at the vertices, to obtain
a CMC surface which is contained in a tubular neighbourhood of this graph.
\item[vi)] Attaching Delaunay cross-pieces \cite{MPPR}: Suppose $\Si$ is a nondegenerate
CMC surface and $p_1,p_2 \in \Si$. Suppose that the tangent planes
at these points are parallel and oppositely oriented (so that the outward
normal of $T_{p_1}\Si$ points towards the outward normal for $T_{p_2}\Si$),
and the distance $|p_1-p_2|$ is an even integer (hence essentially a multiple
of a Delaunay period for $\e$ small). Then, subject to a flexibility condition
on this configuration as well as a slight strengthening of the nondegeneracy condition 
as in iv), one may attach a Delaunay segment along the axis connecting these points.
\end{itemize}

For all of these constructions, with minor provisos concerning the locations
of the points where the gluings are done, the resulting surfaces are nondegenerate,
and hence these operations may be performed iteratively and in various
combinations. I have also omitted a detailed description of the slightly 
strengthened form of nondegeneracy which is needed in iv) and vi). In any case,
using this, it is possible to show, cf.\ \cite{MPP2}, \cite{MPPR}, 
that there exist nondegenerate elements in $\calM_{g,k}$ for any $k \geq 3$
and $g \geq 0$.  

As promised, we sketch the proof of i). Thus let $\Si$ be a nondegenerate 
CMC surface and fix $p \in \Si$.
Excise a small ball $B_\zeta(p)$, and denote by $\Si(\zeta)$ the surface
with boundary, $\Si \setminus B_\zeta(p)$. Now let $D_\e'$ be a half Delaunay
surface, truncated near, but not quite at the neck. More specifically, 
recalling the radial function $\rho_\e(t)$ in the original parametrization of $D_\e$,
which attains its minimum value of $\e$ at $t=0$, we let $D_\e'$ denote
that portion of the Delaunay surface corresponding to $t \geq t_\e$ for 
some small $t_\e < 0$. We now have two surfaces with boundary, and the naive hope
is that if we have truncated them at the correct radii, they should fit
together passably well. This is not quite true without some 
minor modifications. To remedy this, we replace $\Si(\zeta)$ by
a normal graph over it, $\Si(\zeta)' = \Si(\zeta)_{\e G}$, where $G$ is the
Green function for the Jacobi operator with pole at $p$ (this exists
because of the nondegeneracy assumption!). The radius $\zeta$ is chosen
so that the mean curvature of this surface is still bounded. The reason
for doing this is that the shape of the graph of the Green function, hence
of this new surface, is approximately logarithmic, and this matches the
neck region of the Delaunay surface. At any rate, we have now obtained 
two surfaces depending on a parameter which almost match at their boundaries;
for convenience we relabel them as $\Si_1^\e$ and $\Si_2^\e$, respectively
(so the first corresponds to the original surface $\Si$, and the second
to $D_\e^+$).   

The second step of the construction is to consider over each of these surfaces 
separately the space of {\it all} nearby CMC surfaces $(\Si_j^\e)_\phi$ which 
are written as normal perturbations. Here $\phi \in \calC^{2,\al}(\Si_j^\e)$ is
small. This is an infinite dimensional space. For $\Si_1^\e$ it is parametrized 
by arbitrary (small) boundary data $\psi \in \calC^{2,\al}(\del \Si_1^\e)$. 
For $\Si_2^\e$ almost the same is true, except that not every function $\psi$
on the boundary of a half Delaunay surface corresponds to a CMC normal perturbation.
The functions for which this fails lie in the span of the cross-sectional
eigenmodes $e^{i\ell \theta}$, $\ell = 0, \pm 1$. What saves the day is
that these functions are precisely the boundary values of the explicit
geometric Jacobi fields on $\Si_2^\e$, and correspond to CMC surfaces which are rotations
or translations of the original Delaunay surface. Thus we can also regard
all small boundary values on $\del \Si_2^\e$ as corresponding to CMC deformations
of $\Si_2^\e$. The analysis required here is quite simple, and reduces
to a straightforward contraction mapping argument.

The final step is to consider the set of all Cauchy data of these infinite
dimensional spaces of CMC deformations of the summands. The Cauchy data of
a normal perturbation $(\Si_j^\e)_\phi$, by definition, is the 
pair of functions $(\psi_1, \psi_2)$, where $\psi_1$ is the restriction
of $\phi$ to the boundary and $\psi_2 = \del_\nu \phi$ is its normal derivative.
The whole point of this argument is that if we can show that these
Cauchy data subspaces intersect, then that point of intersection corresponds
to CMC deformations of each surface such that these perturbed surfaces match
up to second order along the interface. By elliptic regularity for the mean
curvature equation, these surfaces must actually fit together smoothly,
and we are done. The argument that these infinite dimensional submanifolds
must intersect is also not too difficult. The tangent space of either
submanifold at $\phi = 0$ is the graph of the Dirichlet-to-Neumann
operator for $L_{\Si_j^\e}$. One shows that the
sum of the graphs of these two Dirichlet-to-Neumann operators spans
the whole space when $\e$ is small enough (by verifying
that the same is true for the operators obtained in the limit at $\e=0$).
Finally, the nonlinear Cauchy data submanifolds are graphs off of these
linear subspaces, and some estimates are required to see that the neighbourhood
in which these graphical representations are valid is large enough
to contain the putative point of intersection. 

The proofs of the remaining constructions ii) - vi) may all be done similarly. We 
note that in v) spheres are being used as summands, and these are only nondegenerate
in the broader sense discussed at the very end of \S 2.

\section{Moduli space theory} 
We have defined the spaces $\calM_{g,k}$ as the set of all Alexandrov embedded 
complete CMC surfaces of genus $g$ with $k$ ends. Since we do not mod
out by rigid motions of the ambient space or internal isometries, this might 
reasonably be called the premoduli space, but we shall simply call it the
moduli space of CMC surfaces with given $g$ and $k$. 

We list a few facts about these moduli spaces which are well-known,
or which follow immediately from the material in \S 2 and 3:
\begin{itemize}
\item $\calM_{g,1}$ is empty for all $g \geq 0$, \cite{Me}.
\item $\calM_{g,2}$ is empty when $g \geq 1$. $\calM_{0,2}$ is equal to
the set of all rotations and translations of Delaunay surfaces. These 
follow from Alexandrov reflection arguments, \cite{KKS}.
\item $\calM_{g,k}$ is nonempty for every $g \geq 0$ and $k \geq 3$. 
This was essentially proved by Kapouleas \cite{Ka1}, but also follows from the 
gluing constructions of \S 3. We mention some special cases. Define
$\MM_{g,k} = \{(\Si,p): \Si \in \calM_{g,k},\  p \in \Si\}$. Then there are
continuous mappings
\[
\begin{array}{rcll}
\MM_{g,k} & \longrightarrow&  \calM_{g,k+1} \qquad &\mbox{(from i))}, \\
\MM_{g_1,k_1} \times \MM_{g_2,k_2} & \longrightarrow & \calM_{g_1+g_2,k_1+k_2} \qquad
& \mbox{(from ii))}, \\
\calM_{g,k} & \longrightarrow & \calM_{2g, 2k-2} \qquad &\mbox{(from iv))}.
\end{array}
\]
This final map, which comes from applying the end-to-end construction iv) to two copies
of the same surface, thus doubling it across an end, is only defined on
some open set in the moduli space where some nondegeneracy condition is satisfied.
We mention finally that there are many known minimal $k$-noids of high genus, 
as catalogued in \cite{MPP2},and for each such $g,k$ we obtain by construction iii) an 
element of $\calM_{g,k}$.
\end{itemize} 

The general structure of these moduli spaces is contained in the
\begin{theorem}[\cite{KMP}] For each $g,k$, the moduli space $\calM_{g,k}$ is
a finite dimensional real analytic variety in a neighbourhood of 
each of its points. If $\Si \in \calM_{g,k}$ is nondegenerate,
then there exists a neighbourhood $\Si \in \calU \subset \calM_{g,k}$
which is a real analytic manifold of dimension $3k$.
\end{theorem}
The first assertion in this theorem means that if $\Si \in \calM_{g,k}$,
there is an open neighbourhood $\calV$ of $\Si$ in the space of {\it all}
surfaces near to $\Si$ (in the topology of an appropriate separable Banach
space), and a real analytic diffeomorphism $\Phi: \calV \to \calV'$
such that $\Phi(\calV \cap \calM_{g,k})$ lies in a finite dimensional
linear subspace $B$; furthermore, there exists a real analytic function
$F$ defined on $B \cap \calV'$ such that $\Phi(\calM_{g,k}\cap \calV) = F^{-1}(0)$.
Notice that the dimension of $B$ may depend on $\Si$, and in general 
it is unclear if it is bounded as $\Si$ varies over the moduli space.

When $\Si$ is nondegenerate, the proof is a straightforward application of
Proposition 2 and the implicit function theorem. To prove the more general
assertion, one uses Proposition 1 and the `Kuranishi trick' (Ljapunov-Schmidt 
reduction). 

Note that the generic dimension of this moduli space only depends on the number
of ends, but not on the genus $g$.  

This result illuminates the importance of nondegeneracy in the moduli space
theory. Some consequences of this theorem are that $\calM_{g,k}$ is a real
analytic stratified space, or more precisely, admits a locally finite 
decomposition into real analytic strata. If $\Si$ is a nondegenerate element
in some stratum, then that stratum is maximal, i.e.\ it is not contained in the
closure of any other stratum, and has the predicted dimension $3k$. Unfortunately, this
theorem says nothing about the structure of the moduli space near its
boundary, nor does it limit the number of components it might have. 

We explain the dimension $3k$ which appears here. Supposing that $\Si$
is nondegenerate, the implicit function theorem applies to (\ref{eq:augnl})
since its linearization (\ref{eq:aug}) is surjective. Hence a neighbourhood
of $\Si$ in $\calM_{g,k}$ has a real analytic coordinate chart parametrized
by a ball in the nullspace of (\ref{eq:aug}), and so we must show that
this nullspace is $3k$-dimensional. Although not stated explicitly in \S 2,
this nullspace is precisely the same as the nullspace of $L_\Si$ on
$\calC^{2,\al}_\mu(\Si)$ for small $\mu>0$. This mapping is surjective, by 
nondegeneracy, and the usual integration by parts argument shows that its 
nullspace is identified with the cokernel of $L_\Si$ on $\calC^{2,\al}_{-\mu}$. 
On the other hand, 
\[
6k = \dim W_\Si = \dim \ker \left. L_\Si\right|_{\calC^{2,\al}_\mu} +
\dim \coker \left. L_\Si\right|_{\calC^{2,\al}_{-\mu}}.
\]
These facts together show that the dimension of the nullspace
is $3k$, as claimed. In general, however, and as an alternate derivation
of this formula, we have that
\[
\left. \ker L_\Si\right|_{\calC^{2,\al}_\mu} = \left. \coker L_\Si\right|_{\calC^{2,\al}_{-\mu}},
\qquad 
\left. \coker L_\Si\right|_{\calC^{2,\al}_\mu} = \left. \ker L_\Si\right|_{\calC^{2,\al}_{-\mu}},
\]
where each equality here connotes that the left and right sides are identified by the
integration pairing. Hence
\[
\mbox{ind}\,\big(L_\Si: \calC^{2,\al}_\mu(\Si) \to \calC^{0,\al}_{\mu}(\Si)\big) =
\]
\[
\frac12 \left(\mbox{ind}\,\big(L_\Si: \calC^{2,\al}_\mu(\Si) \to \calC^{0,\al}_{\mu}(\Si)\big)
-  \mbox{ind}\,\big(L_\Si: \calC^{2,\al}_{-\mu}(\Si) \to \calC^{0,\al}_{-\mu}(\Si)\big)\right).
\]
This last quantity can be computed by a relative index theorem \cite{KMP}: the answer
can be computed {\it locally} at infinity. The contribution from each end is $3$, and 
the full relative index is the sum of these contributions over all $k$ ends. 

For any $\Si \in \calM_{g,k}$, Proposition 3 gives a map
\[
\iota: \left. \ker L_\Si\right|_{\calC^{2,\al}_\mu}\longrightarrow W_\Si;
\]
each Jacobi field $\phi$ gets mapped to a $6k$-tuple consisting of the
components of each of the $6$ model geometric Jacobi fields on each of 
the $k$ ends. If $\Si\in \calM_{g,k}$ is nondegenerate, this map is
injective, and hence we can regard $T_\Si \calM_{g,k} = \ker L_\Si$
as a subspace of $W_\Si$. 
\begin{theorem}[\cite{KMP}] There is a natural symplectic structure
on $W_\Si$ such that when $\Si$ is nondegenerate, the $3k$-dimensional
subspace $\iota(\ker L_\Si) \subset W_\Si$ is a Lagrangian subspace.
\end{theorem}

This result is somewhat at odds with the geometric structures
which exist on other standard moduli spaces. In particular, there
does not seem to be a canonical Riemannian metric on these moduli
spaces, and it remains unclear what this Lagrangian structure means,
or how it may be used. Kusner \cite{Ku} has ventured some interesting
speculations about its interpretation. 

At almost the same time these theorems were proved, Perez and Ros \cite{PR} 
discovered the analogous result for certain moduli spaces of complete
minimal surfaces of finite total curvature. Their techniques were 
rooted in the underlying Riemann surface theory, in particular the
Weierstra{\ss} representation.

In the past few years, some new advances have been made concerning these
moduli spaces. The first is a striking result of Gro{\ss}e-Brauckman,
Kusner, and Sullivan.
\begin{theorem}[\cite{GKS1}, \cite{GKS2}] Let $E_3$ be the group of rigid motions of $\RR^3$. Then
$\calM_{0,3}/E_3$ is homeomorphic to the $3$-ball.
\end{theorem}
The argument to prove this relies on a beautiful mix of a classical geometric 
construction relating CMC surfaces in $\RR^3$ with minimal surfaces in $\sphere^3$
and (Dennis Sullivan's) $\ZZ_2$ degree theory for proper mappings in the real analytic category.
This theorem requires the existence of a nondegenerate element in $\calM_{0,3}$, which 
is known from \cite{MP1}. Notice that the dimension count is correct: $\dim E_3 = 6$ 
and so the (pre)moduli space is $6+3=9$ dimensional, which agrees with
the $3 \cdot 3$ dimensions predicted by Theorem 1. 

It is not known whether the homeomorphism in Theorem 3 is a diffeomorphism, nor whether
can exist any triunduloids, i.e.\ elements of $\calM_{0,3}$, which are degenerate.
This latter question is quite important, and if answered in the
negative, would have several interesting implications.

To phrase the final results, we let $\calR_{g,k}$ be the Riemann moduli
space for surfaces of genus $g$ with $k$ punctures, and define the `forgetful map':
\[
\calF: \calM_{g,k} \longrightarrow \calR_{g,k}.
\]
Thus, to any CMC surface $\Si$ we associate its marked conformal completion 
$[\overline{\Si},p_1,\cdots, p_k]$. Thus punctured neighbourhoods around these
points $p_j$ correspond to (slightly perturbed) half Delaunay surfaces.
We also let $\calP_{g,k} = \calR_{g,k} \times \RR^k$ and define an enhanced forgetful map 
\[
\tilde{\calF} : \calM_{g,k} \longrightarrow \calP_{g,k},
\]
\[
\tilde{\calF}(\Si) = \big([\overline{\Si},p_1,\cdots, p_k], \e_1, \ldots, \e_k).
\] 
In other words, this map also records the asymptotic necksizes of each end.
By generalizing the arguments in \cite{KKS}, Kusner has recently proved the
very useful
\begin{theorem}[\cite{Ku}] Up to rigid motions, this enhanced forgetful 
map $\tilde{\calF}$ is proper.
\end{theorem}
This result implies that if $\Si_j$
is any sequence of elements in $\calM_{g,k}$, then 
one of the following three possibilities must hold:
\begin{itemize}
\item There exists a sequence of rigid motions $T_j$ such that $T_j(\Si_j)$
converges to an element $\Si \in \calM_{g,k}$,
\item The marked conformal structures of the sequence $([\overline{\Si}_j],
p_{j,1}, \ldots, p_{j,k})$ must degenerate in $\calR_{g,k}$, or
\item At least one of the necksizes $\e_{j,\ell}$ tends to zero as $j \to \infty$.
\end{itemize}
This gives a very clean and intuitive picture of the modes of possible
degeneration of sequences of CMC surfaces. Kusner also shows that the $\ZZ_2$
degree of this map is zero, so there an even number of preimages of
any regular value.

A topic which warrants much closer examination is the structure of the
moduli spaces $\calM_{g,k}$ near their ends. More specifically,
it would be very interesting to prove that $\calM_{g,k}$ has a
tractable, e.g.\ (semi)analytic, compactification, presumably obtained by
adding moduli spaces $\calM_{g',k'}$ with $g'\leq g$, 
$k' \leq k$, with at least one of the inequalities strict.

It is also natural to try to characterize the image of the unenhanced forgetful
map $\calF$. One motivation for this is that the Teichm\"uller space $\calR_{g,k}$
is known to have fairly complicated topology (its fundamental group is
the Artin braid group on $k$ generators of the topological surface $\overline{\Si}$), 
and so if the image of $\calF$ is relatively large, then $\calM_{g,k}$ itself must have
lots of topology.

\begin{theorem}[\cite{MPP2}] For all $g\geq 0$ and $k\geq 3$,
the mapping $\calF$ is real analytic, and its image in $\calR_{g,k}$ 
is a semianalytic (hence stratified) variety $\calI_{g,k}$. When $g=0$,
then for every $k \geq 3$, $\calF$ is surjective. Likewise, if $g > 0$,
then the codimension of the principal stratum of $\calI_{g,k}$ in $\calR_{g,k}$ 
is uniformly bounded as $k\to \infty$.
\end{theorem}
The main tool in the proof is the half Delaunay addition construction from \S 3;
the proof of the first analyticity assertion reduces to showing that the
uniformizing map which takes $\Si$ to the unique conformally related constant
negative curvature surface is real analytic. The surjectivity statement when $g=0$ 
was also obtained by Kusner \cite{Ku}, and a somewhat weaker statement can 
also be deduced from Ratzkin's end-to-end gluing construction \cite{Rat}. 

There are many interesting questions raised by these theorems.
Kusner \cite{Ku} lists several important ones concerning the finer
structure of the map $\tilde{\calF}$, and in particular concerning
the precise number of points in the preimages of regular values.   

\section{CMC surfaces in other asymptotically Euclidean $3$-manifolds}
Let $(X,h)$ be a complete three dimensional manifold with asymptotically
Euclidean ends. Recall that this means that the following hypotheses hold:
there exists a compact set $K \subset X$ such that $X \setminus K$ is a disjoint 
union of ends, $F_1, \ldots, F_\ell$, each of which is equipped with a diffeomorphism
$\Psi_j:\RR^3 \setminus B(0,R_j) \to F_j$ such that the metric $h$ on $X$
pulls back via $\Psi_j^*$ to a metric on the exterior region in $\RR^3$ of the form
$\delta + {\mathcal O}(r^{-1})$, with corresponding decay for its derivatives. 

We shall consider the space $\calM_{g,k}(X,h)$ consisting of all proper
Alexandrov embedded CMC surfaces in $X$ of genus $g$ with $k$ ends. If $X$ has more than one end,
then we could also consider a decomposition into moduli spaces which
explicitly label which end of $\Si$ tends to infinity in which end of $X$, but we
shall not introduce special notation for this. As noted in \S 2, the asymptotics 
result of \cite{KKS} most likely extends to this setting, but since this
is not worked out anywhere, we circumvent this by assuming that this moduli space 
contains only CMC surfaces which are asymptotically Delaunay. 

It is straightforward to check that the main results about the structure of the 
moduli space for CMC surfaces in $\RR^3$ persist in this more general situation. 
In particular, we still have that for any $(X,h)$, $\calM_{g,k}(X,h)$ is a locally 
real analytic variety, and near nondegenerate points it has dimension $3k$. Note 
that we do not need to assume that $X$ or $h$ are real analytic because this theorem 
relies on the implicit function theorem in some function space on $X$, and the mean 
curvature functional depends real analytically on the metric and its derivatives.

Consider first the case where $X = \RR^3$ and $h$ is a small (decaying) perturbation of the 
standard metric $\delta$. We prove below that $\calM_{g,k}(\RR^3,h)$ is a small 
perturbation of the `standard' moduli space $\calM_{g,k}(\RR^3,\delta)$, as well as
a standard 'generic regularity' theorem, that this perturbed moduli space is 
smooth for generic $h$; in other words, for most small alterations of the ambient metric  
there are no degenerate CMC surfaces. 

As an aside, I have already noted that no degenerate CMC surfaces are known
to exist, and it seems very unlikely that it will be possible to construct
one directly. If they exist at all, I suspect they will only be found by
indirect means, by showing that the projection mapping $\Pi$ of the
big moduli space $\calS$ defined below must have fold singularities.
It is completely unclear whether this happens for the standard (or any!)
metric on $\RR^3$. 

Let $Z \equiv S^{2,\al}_{-1}(\RR^3)$ denote the space of all symmetric $2$-tensors with 
coefficients in the weighted H\"older space $\calC^{2,\al}_{-1}(\RR^3)$. Thus if $x$
is the standard Cartesian coordinate chart on $\RR^3$, then any $\gamma \in Z$ satisfies 
$|\del_x^\beta \gamma_{ij}| \leq C |x|^{-1-|\beta|}$ as $|x| \to \infty$,  $|\beta| \leq 2$, 
and with the H\"older quotients of the second derivatives decaying like $|x|^{-3-\alpha}$. 
Let $\calU$ be an open ball in this space such that if $\gamma \in \calU$, then $\delta + \gamma$ is 
everywhere positive definite, hence is an asymptotically Euclidean metric on $\RR^3$. 

\begin{theorem} There exists a dense subset $\calU' \subset \calU$ such that when 
$\gamma \in \calU'$, $\calM_{g,k}(\RR^3,\delta+\gamma)$ is a smooth analytic manifold of dimension $3k$.
\end{theorem}
\begin{proof} Consider the Banach manifold $\calB$ of all Alexandrov embedded $\calC^{2,\al}$ 
surfaces in $\RR^3$ of genus $g$ with $k$ asymptotically Delaunay ends. A coordinate chart 
on this manifold near any point $\Si$ is given by a small ball in $W_\Si 
\oplus \calC^{2,\al}_{-\mu}(\Si)$. Thus $w \in W_\Si$ induces a change of the ends, 
either by rigid motion or change of necksize, and $u \in \calC^{2,\al}_{-\mu}(\Si)$ 
gives a normal perturbation. We write the resulting surfaces as $\Si_{w,u}$. 

The first step is to show that the set
\[
\calS = \{(\Si,\gamma) \in \calB \times Z: \Si \in \calM_{g,k}(\RR^3,\delta + \gamma)\}
\]
has the structure of a Banach submanifold near $\gamma=0$. This is a straightforward 
application of the implicit function theorem. Fix $\Si \in \calM_{g,k}$ and a
coordinate chart $\calV$ on $\calB$, and define the functional 
$N: \calV \times Z \to \calC^{0,\al}_{-\mu}(\Si)$
by setting $N(w,u,\gamma)$ equal to the mean curvature function on $\Si_{w,u}$ with respect to 
the metric $\delta + \gamma$. The differential of this mapping $D_{12}N$ in the first
two slots, at $(w,u)=(0,0)$, is the Jacobi operator $L_\Si$, and by Proposition 2 this
operator has closed range of finite codimension. We must show that the full differential 
$DN$ at $(0,0,0)$ is actually surjective. Granting this for the moment, we then obtain 
an analytic mapping $\Psi$ defined in a neighbourhood of $0$ in $Z$ to $\calV$ such that 
all points in $\calS$ near to $\Si$ are of the form $(\Psi(\gamma),\gamma)$. 

To prove the claim about surjectivity, we must show that 
\[
\mbox{ran}\,\big(L_\Si\big) + \mbox{ran}\,\big(\left.D_3N\right|_{0}\big) 
= \calC^{0,\al}_{-\mu}(\Si).
\]
It is necessary to reinterpret this slightly. First, the mapping (\ref{eq:aug})
is the restriction of the mapping (\ref{eq:mL}), i.e.\ from 
$\calC^{2,\al}_{\mu}(\Si)$ to $W_\Si \oplus \calC^{2,\al}_{-\mu}(\Si)$. 
By Proposition 3, the nullspaces of these mappings are the same, hence their
cokernels also agree. As in \S 4, duality considerations show that 
\[
\calC^{0,\al}_{\mu} = \mbox{ran}\,\big(\left. L_\Si\right|_{\calC^{2,\al}_{\mu}}\big)
\oplus \ker \left. L_\Si \right|_{\calC^{2,\al}_{-\mu}}.
\]
Hence we must show that there does not exist any $\psi \in \calC^{2,\al}_{-\mu}(\Si)$
such that $L_\Si \psi = 0$ and $\langle \psi , D_3N(\gamma) \rangle = 0$
for all $\gamma \in Z$. In the end, we show that when $\gamma$ ranges over all $\calC^{2,\al}$
symmetric $2$-tensors supported in some small ball $\calW$ in $\Si$, the
functions $D_3N(\gamma)$ fill out $\calC^{0,\al}(\calW)$, and this will prove
the claim since if the orthogonality condition were to hold, then $\psi$
would have to vanish in this set $\calW$, which is impossible since it
is in the nullspace of $L_\Si$. 

By dilating the space substantially and using an approximation argument, 
it will be enough to consider the simpler problem where $\Si = \RR^2 \subset
\RR^3$, and so we must compute the linearization $D_3 N$ which measures
how the mean curvature of the flat plane $\{(x_1,x_2,0)\}$ in $\RR^3$ changes as we vary
the metric. This computation is not too horrible, fortunately, and we
obtain that
\[
D_3 N(\gamma) = (\mbox{div}\, \gamma )_3 + \frac{\del\,}{\del x_3}\frac12 (\mbox{tr}\, \gamma);
\]
here $\gamma = (\gamma_{ij})$ and the divergence and trace are computed only
in the $(x_1,x_2)$ directions. The subscript $3$ in the first term on
the right denotes the third component of $\mbox{div}\,\gamma$, i.e.\ $\gamma_{i3;}^{\ \ \ i}$.
Examining this expression, it is clear that we can specify it 
arbitrarily in the bounded set $\calW$. This proves the claim. 

Consider the projections $\Pi_1: \calB \times Z \to \calB$ and $\Pi_2: \calB \times Z \to Z$.
The restriction of $\Pi_2$ to $\calS$ has surjective differential at $\Si'$ if and only $\Si'$
is nondegenerate with respect to the metric $\delta + \gamma$, $\gamma = \Pi_2(\Si')$. 
Hence by the Sard-Smale theorem, if $K\subset\calB$ is a compact set containing $\Si$,
then there is an open dense set $\calU_K \subset Z$ of metrics near $\delta$ such that
every surface in $\Pi_1^{-1}(K) \cap \calS \cap \Pi_2^{-1}(\calU_K)$ is nondegenerate
(with respect to the appropriate metric).

To globalize this, cover $\calM_{g,k}$ by a countable union of compact 
sets $\calK_j$ and choose open dense neighbourhoods $\calU_j$ in $Z$ as above.
Then the intersection of these sets $\calU_j$ is still dense, by Baire's 
theorem, and has the desired property. 
\end{proof}

It is most likely possible to modify this proof to show something 
stronger and much more useful. Let $\calD$ denote the set of `bad' metrics, 
i.e.\  those for which the perturbed moduli space contains degenerate elements, 
then I expect it is true that $\calD$ decomposes into a union of $\calD' \cup \calD''$ 
where $\calD'$ is a codimension $2k$ submanifold in $\calU$ and 
$\calD''$ is quantitatively smaller. (In fact, probably $\calD$ has 
a stratified structure with only strata of bounded codimension.) The reason 
for this conjecture is as follows. We may view this process of destroying degenerate
surfaces as a problem in eigenvalue perturbation theory. As explained in
\S 2, the essential spectrum of $L_\Si$ is a locally finite union of intervals
$I_j$ and $0$ always lies on the left endpoint of one of these intervals.
In addition, there is always a small interval $(-\e,0)$ which is disjoint from
the spectrum; here $\e$ depends on $\Si$. $\Si$ is degenerate if and only if 
$0$ is also in the point spectrum. 
There is a good procedure for tracking what happens to this $L^2$ eigenvalue as 
$\Si'$ varies in $\calB$. To explain it, consider (for any fixed $\Si'$) the resolvent 
$R(\lambda) = (-L_{\Si'} - \lambda^2)^{-1}$; if we restrict $\lambda^2$ to lie in 
$B(0,\e)\setminus [0,\e)$, or equivalently, if $\lambda$ lies in the half-ball 
$B(0,\sqrt{\e}) \cap \{\mbox{Im} \lambda < 0\}$, which is contained in the resolvent
set of $-L_{\Si'}$, then $R(\lambda)$ is bounded on $L^2(\Si')$. Using 
the machinery developed in \cite{MPU} it is quite straightforward to show
that $R(\lambda)$ continues meromorphically to some neighbourhood of $0$
in either the complex plane or the logarithmic complex plane (the latter 
being necessary if $0$ is a nontrivial ramification point in this continuation). 
This continued resolvent has at most a double pole at the origin, and the
coefficient of $1/\lambda^2$ is the projection onto the space of $L^2$
Jacobi fields. Thus we are interested in the dependence of this meromorphic
structure as $\Si'$ varies. If we were starting with some eigenvalue $\lambda_0^2 > 0$
(i.e.\ inside the continuous spectrum, but away from the `threshold value' $0$), then 
$R(\lambda)$ has only a simple pole at $\lambda_0$. As we deform $\Si$, this pole could move
either to the left or right on the real axis or upwards, into the so-called
nonphysical half-plane (where it would become a resonance, or scattering pole).
This is predicted by `Fermi's golden rule' \cite{RS}. Agmon and Herbst 
have studied persistence of embedded eigenvalues and have shown in 
as yet unpublished work, but cf.\ also \cite{CHM}, that for non-threshold eigenvalues, 
the set of perturbations which leave the eigenvalue fixed is a submanifold of 
codimension equal to the multiplicity 
of the continuous spectrum, which in this case equals $2k$. Unfortunately,
their proof does not carry over in an obvious way to our case, where
the pole is double and occurs at the ramification point. I hope that 
this will be possible, and if it is, full details will be given elsewhere.

This somewhat arcane issue is interesting for the following reason. 
We have shown that for most $h$ near $\delta$, $\calM_{g,k}(\RR^3,h)$
is a smooth (in fact, analytic) manifold of dimension $3k$. However,
the topology of this manifold may depend on $h$. It would be quite
interesting if there were a canonical smoothing of the moduli
space $\calM_{g,k}(\RR^3,\delta)$, and to show this it would be sufficient to 
know that the set $\calD$ of bad metrics has codimension at least two, so that
there are no `wall-crossing' transitions. This 
is precisely the point of the preceding discussion.  We note that this discussion
carries over to the more general setting of the moduli space of CMC 
surfaces in general asymptotically Euclidean manifolds $(X,h)$, and
so we might then expect to get canonical smooth moduli spaces which depend 
only on $X$, but not on $h$!  A case of particular interest is when
$X = \overline{X} \setminus \{p_1,\ldots, p_\ell\}$, the complement of an
ordered finite set of points in a smooth compact manifold, endowed with
a metric which is asymptotically Euclidean near these punctures. 
We would then obtain a canonical  
object $\calM_{g,k}(\overline{X},\{p_1,\ldots, p_\ell\})$ associated only 
to the smooth structure of $\overline{X}$ and the isotopy class of the 
marked set of points.  As ventured earlier, one might hope that this CMC 
moduli space could reflect some of the topology of this data. 

We conclude by discussing a final example: if $p \in \overline{X}$, 
and $h$ is an asymptotically Euclidean metric on $X = \overline{X} \setminus \{p\}$,
then $\calM_{0,2}(X,\e^{-2}h)$ contains many nontrivial elements when $\e$ is
small enough. The construction is quite simple. Let $\gamma$ be a bi-infinite 
geodesic which is properly embedded in $X$; thus both ends converge to $\infty$
in the punctured ball around $p$. Choose Fermi coordinates $(s,r,\theta) \in
\RR \times (0,1) \times S^1$ around $\gamma$ corresponding to the metric $\e^{-2}h$ 
and use these to define the approximate Delaunay surface $r = \rho_\e(s)$.
The mean curvature of this is approximately $1$, but is not exactly constant.
A straightforward perturbation argument shows that there is a small perturbation of
this surface which is CMC; we omit the details. Notice that there may
be many nonisotopic geodesics, depending on the topology of $X$, and hence
there may be many components in this moduli space. In fact, one expects that
the components are in one-to-one correspondence with $\pi_1(X,p)$.

Once these Delaunay-type surfaces have been constructed, it is also not
difficult to carry out the analogue of the construction i), of adding
Delaunay ends, and most of the other constructions too. 

While it is far from clear that this extension of the CMC moduli space 
theory will have any uses, these spaces are certainly natural geometric objects
and I suspect that there is a lot more structure which has yet to be
uncovered.

\end{document}